\documentclass[preprint,12pt]{elsarticle}

\usepackage{amsmath,amssymb}

\newtheorem{thm}{Theorem}[section]

\newtheorem{prob}[thm]{Problem}

\journal{Expositiones Mathematicae}

\begin{document}
	
	\begin{frontmatter}
		\title{A historical perspective of Tian's evolution algebras}
		
		\author[1]{Manuel Ceballos}
		\ead{mceballos@uloyola.es}
		\author[2]{Ra\'ul M. Falc\'on\corref{cor}}
		\cortext[cor]{Corresponding author.}
		\ead{rafalgan@us.es}
		\author[3]{Juan
			N\'u\~nez-Vald\'es}
		\ead{jnvaldes@us.es}
		\author[4]{\'Angel F. Tenorio}
		\ead{aftenorio@upo.es}

		\address[1]{Departamento de Ingenier\'ia. Universidad Loyola Andaluc\'ia. Seville, Spain.}
		\address[2]{Departmento de Matem\'atica Aplicada I, Universidad de Sevilla. Seville, Spain.}
		\address[3]{Departamento de Geometr\'{\i}a y Topolog\'{\i}a. Universidad
			de Sevilla. Seville, Spain.}
		\address[4]{Dpto. de Econom\'ia, M\'etodos Cuantitativos e Historia
			Econ\'omica. Universidad Pablo de Olavide, Seville, Spain.}

		\begin{abstract} Even if it has been less than a decade and a half since Tian introduced his concept of evolution algebras to represent algebraically non-Mendelian rules in Genetics, their study is becoming increasingly widespread mainly due to their applications to many scientific disciplines. In order to facilitate further research on the topic, this paper deals with the past and present research on these kind of algebras, together with the most
			relevant topics regarding them.
		\end{abstract}
		
		
	\end{frontmatter}
	
	%\linenumbers

	\section{Introduction}
	
	Dealing with the chloroplast inheritance in plants, the German botanists and geneticists Carl Correns (1864--1933) \cite{Correns} \cite{Correns1909} and Erwin Baur (1875--1933) \cite{Baur1909} realized the relevance that non-Mendelian rules have in Genetics. In spite of this, and unlike Mendelian Genetics, for which an algebraic interpretation of its rules was already introduced in the 1930's by Serebrowski \cite{Sere} and Kostitzin \cite{Kos}, it is not until the early 2000s that a similar interpretation for non-Mendelian rules was proposed.
	
	More specifically, in order to relate Markov chains with the prokaryotic cell reproduction, and based on the self-reproduction rules of non-Mendelian Genetics, evolution algebras were firstly introduced in the Ph.D. Thesis of Jianjun Paul Tian \cite{tesistian} in 2004; then jointly presented with Petr Vojtechovsky \cite{Petr} in 2006, and later analyzed in greater depth in a book by Tian \cite{Tian2008} in 2008. The dynamic nature of these non-associative algebras is not described by a series of identities, unlike many other known non-associative algebras such as alternative, Lie, Malcev or Jordan algebras, amongst others. More specifically, an {\em evolution algebra} $E$ over a field $\mathbb{K}$ is a finite-dimensional algebra such that there exists a {\em natural basis} $\{e_1, \ldots, e_n\}$ so that $e_i \cdot  e_j = 0$, for all $i,j\in\{1,\ldots,n\}$ such that $i \neq j$, and $e_i \cdot  e_i = \sum_k \, p_{ik} \, e_k,$ for all $i\in\{1,\ldots,n\}$, and some coefficients $p_{ik} \in \mathbb{K}$, which are called the {\em structural constants} of the evolution algebra $E$. This multiplication is extended linearly from the given multiplication of basis elements. By considering each generator $e_i$ as an allele in Genetics, the {\em structure matrix} $(p_{ij})_{i,j}$ encodes somehow the laws of inheritance of non-Mendelian Genetics as algebraic properties of the algebra. 
	
	Even if a wide range of mathematicians has dealt with this type of algebras since the original manuscripts of Tian and Vojtechovsky, the relevance of evolution algebras is not still widely known. In order to make easier further research on the topic, this paper gathers together the antecedents, origin, early stage and later development of this kind of algebras, and shows different applications arising from them. In any case, we recommend the recent manuscript of Utkir Rozikov \cite{Rozikov2020} for some detailed descriptions on genetic algebras, evolution algebras and quadratic stochastic operators, with particular emphasis in the study of algebraic and probabilistic approaches in population dynamics. See also \cite{Qaralleh2017} for a short survey on genetic and evolution algebras, and the expository article writen by Tian \cite{invita} in 2016.
	
	The paper is organized as follows. Section \ref{sec:Antecendents} outlines the origin and development of genetic algebras as primordial antecedents of Tian's concept of evolution algebras. In Section \ref{sec:Preliminaries} some preliminary concepts and results on evolution algebras are indicated.  Next, in each subsection of Section \ref{sec:FirstSteps}, a brief comment on the main papers published on this topic is given. Furthermore, Section \ref{sec:Topics} overviews the relationship among evolution algebras, Graph theory, Group theory, Markov chains and Biology. The paper is completed with an extensive bibliography, which may be of valuable help to all those researchers interested in this topic. Throughout the paper, we follow the current notation, which may differ from the original one to which we refer.
	
	\section{Antecedents} \label{sec:Antecendents}
	
	This section deals with the origin of genetic algebras as cornerstone of Tian's concept of evolution algebra. For more details about the historical background on genetic algebras, we refer the reader to \cite{Bert1966, Lyu1992, Reed1997, Rozikov2020, Worz1980}.
	
	\vspace{0.2cm}
	
	In 1866, Gregor Johann Mendel (1822--1884) was the first who made use of algebraic symbols to express sexual reproduction laws of inheritance. It was done in his original manuscript on experiments in plant hybridization \cite{Men}. More specifically, given a pair of plant characters ${\bf A}$ and ${\bf a}$, together with the hybrid form ${\bf Aa}$ in which both characters fuse together, Mendel indicated with an expression of the form ${\bf \alpha A+\beta Aa + \gamma a}$ the ratio $\alpha:\beta:\gamma$ of numbers of offsprings having each possible plant character  ${\bf A}$, ${\bf Aa}$ and ${\bf a}$ in the following generation. In this way, Mendel introduced three principles:
	\begin{itemize}
		\item The {\em principle of paired factors}, which states that characters are controlled by unit factors that appear in pairs within each individual organism.
		\item The {\em principle of dominance}, which states that one of the two previous factors is dominant and makes its effect in the individual, whereas the other one is recessive and does not show its effect unless both characters are recessive.
		\item The {\em principle of segregation}, which states that, during the inheritance process, the mentioned pair of factors separate randomly so that the offspring receives exactly one factor from each parent.
	\end{itemize}
	
	\vspace{0.2cm}
	
	Almost forgotten for decades, Mendel's laws were independently rediscovered around 1900 by the Dutch botanist Hugo de Vries (1848--1935) \cite{Vries05}, who introduced the term {\em mutation} and suggested that of {\em gene}; the German botanist and geneticist Carl Correns (1864--1933) \cite{Correns}, who discovered the cell extranuclear inheritance; the Austrian agronomist Erich von Tschermak (1871--1962) \cite{Tschermak}, who introduced the combination of plant characters in order to improve the efficiency of plant breeding; and the American agronomist William Jasper Spillman (1863--1931) \cite{Spillman09}, whose work was crucial for the development of {\em Agricultural Economics}. The emergence of all their contributions concerning the distribution of characters among offspring in plant hybridization made fundamental the study and development of Mendel's laws.
	
	In 1905, the English biologist William Bateson (1861--1926) introduced the term {\em Genetics} to describe the study of inheritance processes \cite{Bateson}. Bateson was a fervent advocate of Mendel \cite{Bateson02}, whose principles, even being reborn, were questioned at that time. An interesting question in this regard had been asked in 1902 by the British statistician George Udny Yule (1871--1951), who contemplated \cite{Yule02} that dominant factors could increase indefinitely in a continuous way throughout generations. One year later, the English mathematician and biostatistician Karl Pearson (1857--1936) \cite{Pearson03} and the American zoologist and geneticist William Ernest Castle (1867--1962) \cite{Castle03} established that, under random circumstances, there exists certain stability of characters in the inheritance process. Concerning such a stability, the English mathematician Godfrey Harold Hardy (1877--1947) \cite{Hardy08} asked in 1908 about the circumstances under which the distribution of characters in the offspring is the same that in the generation before in the absence of disturbing factors.
	
	The answer was independently given at that moment by Hardy himself and the German physician Wilhelm Weinberg (1862--1937) \cite{Weinberg08}, and is currently known as the {\em Hardy-Weinberg law}. An important limitation of that law is the assumption of a potentially infinite population, not taking into account the importance of sampling fluctuations in evolutionary processes described by finite populations. The first proposal dealing with such a possibility would be the so-called {\em Wright-Fisher model} based on the original manuscripts of the American geneticist Sewall Green Wright (1889--1988) \cite{Wright1931} and the British statistician and geneticist Ronald Aylmer Fisher (1890--1962) \cite{Fisher1922}.
	
	In 1923-24, the Soviet mathematician Sergei Natanovich Bernstein (1880-1968) \cite{Ber,Ber24} introduced the concept of {\em quadratic stochastic operator} (QSO) as a map $V:S\rightarrow S$, where
	\[S:=\{(x_1,\ldots,x_n)\in\mathbb{R}^n\colon\, x_i>0, \text{ for all } i, \text { and } \sum_ix_i=1\},\] 
	and, for each $(x_1,\ldots,x_n)\in S$, it is $V((x_1,\ldots,x_n))=(V(x_1),\ldots,V(x_n))$, where
	\[V(x_k):=\sum_{i,j}p_{ij}^k x_ix_j, \text{ for all } k\in\{1,\ldots,n\}.\]
	Here, the product $x_ix_j$ is the usual product of real numbers. In addition, $p_{ij}^k\geq 0$ and $p_{ij}^k=p_{ji}^k$, for all $i,j,k\in\{1,\ldots,n\}$, and $\sum_k p_{ij}^k =1$, for all $i,j\in\{1,\ldots,n\}$. As such, every QSO constitutes an {\em evolutionary operator} that describes the time evolution or inheritance process of a {\em free population} with $n$ different genetic types.
	
	More specifically, each component $x_i$ of a given element $x=(x_1,\ldots,x_n)\in S$ represents the probability that a random individual in the population under consideration belongs to the species that is determined by the $i^{\mathrm{th}}$ genetic type. Hence, the $n$-tuple $x$ describes the distribution of the population with respect to the $n$ genetic types, whereas $V(x)$ describes such a distribution for the next generation. In particular, each value $p_{ij}^k=p_{ji}^k$ determines the probability that an offspring with genetic type $k$ arises from two individuals of respective genetic types $i$ and $j$, without sexual differentiation.
	
	Based on the possible fluctuation of genetic types throughout subsequent generations of an evolutionary process, a main problem in the theory of QSOs consists of determining their limit behavior for any given initial distribution of genetic types $x\in S$. That is, the study of the corresponding distribution $V^m(x)$ for the $m^{\mathrm{th}}$ generation, when $m$ tends to infinity. In this regard, Berstein focused in particular on the problem of determining and classifying all those QSOs for which a stationary distribution or stable evolution arises in only one generation (that is, such that $V^2=V$). Even if Berstein only solved this problem for $n=3$, it is currently solved for all dimensions (see \cite{Gonzalez95, Gutierrez00}). To this end, it would be crucial the work in 1975 of Philip Holgate (1934--1993) \cite{Holgate75}, who expressed algebraically this problem as follows.
	\begin{enumerate}
		\item Firstly, Holgate related each given evolutionary operator $V$ with the algebra described as
		\[xy=\frac 12 \left(V(x+y)-V(x)-V(y)\right).\]
		(It was called {\em evolution algebra} by Joseph Bayara \cite{bayara99, bayara}, which differs from the concept introduced by Tian, the one that is commonly used in the current literature and with which this paper deals.)
		
		\item Then, he introduced the concept of {\em Bernstein algebra} as an algebra $A$ such that $x^2 \, x^2 = \omega^2(x) \, x^2$, for all $x \in A$, where $\omega$ is an algebra homomorphism from $A$ into its base field, which is called the {\em weight} of the algebra. Recall that an algebra for which one such a non-trivial homomorphism exists is called {\em baric} \cite{Ethe}.
	\end{enumerate}
	
	In 1934, Alexander Pawlowitsch Serebrowski (1884-1938) \cite{Sere}
	was the first to give an algebraic interpretation of the symbol
	$\times$ as a mathematical way to represent Mendelian inheritance
	laws. A similar symbolic multiplication was independently
	introduced shortly after by Aleksandrovich Kostitzin (1883-1963)
	\cite{Kos}. 
	
	At the same period, Valery Ivanovich Glivenkov
	(1896-1940) \cite{Gli} introduced the so-called {\em Mendelian
		algebras} for diploid species. Nevertheless, it was Ivor Malcolm
	Haddon Etherington (1908-1994) who, in 1939, introduced in
	Genetics \cite{Ethe} the systematic study of commutative
	non-associative linear algebras, by describing in particular the
	following three types of algebras:
	\begin{itemize}
		\item A {\em gametic algebra} is a finite-dimensional algebra of basis $\{e_1,\ldots,e_n\}$ such that $e_i\cdot e_j=\sum_kp^k_{ij}e_k$, where each structural constant $p^k_{ij}$ belongs to the real interval $[0,1]$ so that $\sum_{k}p_{ij}^k=1$, for all $i,j\in\{1,\ldots,n\}$. As such, each value $p^k_{ij}$ represents the probability that an arbitrary gamete of zygotic type $e_k$ derives from an individual of zygotic type $e_ie_j$. As occurs with the QSOs, the probability or progeny distribution describes the evolution of the population under consideration. The latter constitutes a free population if the inheritance process  derives from random matings; that is, with absence of sexual differentiation and selection.
		\item A {\em zygotic algebra} is the {\em duplicate} of a gametic algebra. That is, the former is isomorphic to the set of quadratic forms of the latter.
		\item A {\em copular algebra} is the duplicate of a zygotic algebra.
	\end{itemize}
	Each one of the three just described algebras are baric by means of the weight $\omega$ that is linearly defined from $\omega(e_i) = 1$, for all $i\in \{1,\ldots,n\}$. Here, $\{e_1,\ldots,e_n\}$ is the basis of the algebra under consideration, and hence, 
	\[\omega(e_i)\omega(e_j)=1=\sum_k p_{ij}^k=\sum_k p_{ij}^k\omega(e_k)=\omega\left(\sum_k p_{ij}^k e_k\right)=\omega(e_ie_j),\]
	for all $i,j$. Further, Etherington called {\em train algebra} any baric algebra $A$ with weight $\omega$ such that the equation of lowest degree relating the principal powers of any vector $x\in A$ has the form
	\[x^m + c_1 \,\omega(x) \,x^{m-1} + \ldots +  c_m \,\omega(x)^m = 0.\]
	Finally, he termed {\em special train algebra} any baric algebra $A$ with weight $\omega$ such that $N = \mathrm{ker}(\omega)$ is nilpotent (that is, there exists a positive integer $m\in\mathbb{N}$ such that $N^m=0$) and all the principal powers $N^i,$ with $i \in \mathbb{N}$, are ideals of $A$. In particular, every special train algebra is a train algebra.
	
	Etherington called {\em genetic algebra} any one of the family of gametic, zygotic and copular algebras and ensured that {\em ``all the  fundamental genetic  algebras are special  train  algebras''}. Shortly after, he indicated \cite{Ethe4} that this statement is only valid for gametic algebras, but not even for zygotic algebras, because the  duplicate of a special train algebra, although a train algebra, is not always a special train algebra.  This fact provided serious problems for translating some properties and relations from genetic inheritance processes.
	
	In order to avoid the flaws in Etherington's genetic algebras and give rise to a transparent structure theory for them, Richard Donald Schafer (1918--2014) \cite{Schafer} provided in 1949 an alternative formal definition that is intermediate between the concepts of train algebra and special train algebra. To this end, he made use of the {\em transformation algebra} $T(A)$ of a non-associative algebra $A$. It is the algebra of all the polynomials with coefficients in the base field, whose variables are transformations in the set formed by the identity map in $A$, all the {\em right multiplications} $R_{\alpha}: A \mapsto A$ such that $R_{\alpha}(x) = x\, \alpha$,  and all the {\em left multiplications}
	$L_{\alpha}: A \mapsto A$ such that $L_{\alpha}(x)=\alpha\, x$, for all $\alpha, x \in A$. In particular, the transformation algebra of a baric algebra is baric.
	
	Further, unlike Etherington, Schafer assumed the commutativity that is inherent in any inheritance process and hence, he considered $R_{\alpha}=L_{\alpha}$, for all $\alpha\in A$. Then, he called {\em genetic algebra} any commutative baric algebra $A$ such that all the coefficients of the characteristic function $|\lambda I - T|$, with $T\in T(A)$, only depend on the weight of the algebra. As such, every genetic algebra is a train algebra and every special train algebra is a genetic algebra. Moreover, unlike special train algebras, the duplicate of a genetic algebra is always a genetic algebra.
	
	In 1971, Harry Gonsh\"or \cite{Gon2} proved that the concept of genetic algebra proposed by Schafer is equivalent to have a commutative algebra $A$ of basis $\{e_0,e_1,
	\ldots, e_n\}$ such that $e_i \, e_j = \sum_k \, c_{ij}^k \ e_k$, for all $i,j\in\{0,\ldots,n\}$, where $c_{00}^0 = 1$; $c_{0j}^k = 0,$ for all $k<j$; and $c_{ij}^k = 0,$ for all $i,j>0$ and $k\leq \max\, \{i,j\}$. This new definition was proved to have significant relevance in Genetics thanks to the genetic meaning that Holgate \cite{Holgate1972, Holgate1987} gave of duplicates and derivations of a genetic algebra. The study of the subspace $\mathrm{Der}(A)$ of derivations of a genetic algebra $A$ has also been dealt with by Gonsh\"or \cite{Gon1988} himself and also by Roberto Costa \cite{Costa1982, Costa1983}, Aribano Micali and Philippe Revoy \cite{Micali1986}, and, much more recently, by Rasul Ganikhodzhaev, Farrukh Mukhamedov, Abror Pirnapasov and Izzat Qaralleh \cite{Ganikhodzhaev2018}. Recall in this regard that $D\in \mathrm{Der}(A)$ if and only if $D(xy)=D(x)y+xD(y)$, for all $x,y\in A$. As for any algebra, such a subspace constitutes a Lie algebra.
	
	Holgate also introduced \cite{Holgate1} both notions of {\em sex differentiation algebra} and {\em dibaric algebra}. The former characterizes the equal division of offspring between both sexes. More specifically, the sex differentiation algebra is a bi-dimensional commutative algebra of basis $\{m,w\}$ such that $m^2=w^2=0$ and $mw=wm=(m+w)/2$. Further, an algebra is called dibaric if it admits a homomorphism onto the sex differentiation algebra. Holgate proved in particular that, if $A$ is a dibaric algebra, then its {\em derived algebra} $A^2$ is a baric algebra. With the introduction of dibaric algebras, he formalized the original idea of Etherington \cite{Ethe} of treating in a separate way male and female components of the population. 
	
	Another contribution of Holgate \cite{Holgate66} in the theory of genetic algebras was the use of isotopisms of algebras in order to represent algebraically the mutation of genotypes in the inheritance process. Recall in this regard that two $n$-dimensional algebras $A$ and $A'$ are said to be {\em isotopic} \cite{Albert42} if there exist three non-singular linear transformations $f$, $g$ and $h$ from $A$ to $A'$ such that $f(u)g(v) = h(uv)$, for all $u,v\in A$. The triple $(f,g,h)$ is called an {\em isotopism} (an {\em isomorphism}, if $f=g=h$; and an {\em homomorphism}, if besides, the condition of non-singularity is not imposed) between the algebras $A$ and $A'$ (see \cite{FFN18} for a recent survey on the theory of isotopisms). Together with Tania Campos \cite{Campos87}, Holgate proved that certain types of zygotic algebras representing chromosome segregation and recombination are isotopic.
	
	Much more recently, Manuel Ladra, Bakhrom Omirov and Rozikov \cite{Ladra2014} would introduce the concept of {$bq$-homomorphism} as a pair of linear maps $f$ and $g$ from the algebra to the set of real numbers, so that $f(uv)=g(uv)=(f(u)g(v)+f(v)g(u))/2$, for all $u,v\in A$. (It constitutes an homomorphism whenever $f=g$.) Then, they proved \cite[Theorem 3.7]{Ladra2014} that an algebra is dibaric if and only if it admits a non-zero bq-homomorphism.
	
	\section{Preliminaries on evolution algebras}\label{sec:Preliminaries}
	
	In order to make this article as self-contained as possible and facilitate a better understanding for the reader, we recall in this section some notions and basic results on evolution algebras that were introduced in the original manuscripts of Tian and Vojtechovsky \cite{Tian2008, Petr}. To this end, and from now on, let $E$ denote an $n$-dimensional evolution algebra over a base field $\mathbb{K}$, with natural basis $\{e_1,\ldots,e_n\}$ and structural constants $p_{ik}\in\mathbb{K}$, for all $i,k\in\{1,\ldots,n\}$.
	
	If the base field $\mathbb{K}$ is the real field, then the evolution algebra is said to be {\em real}. In such a case, it is said to be {\em non-negative} if $p_{ik}\geq 0$, for all $i,k\in\{1,\ldots,n\}$. If, besides,  $\sum_k p_{ik}=1$, for all $i$, then the real and non-negative evolution algebra is called {\em Markov evolution algebra}.
	
	Every evolution algebra $E$ is non-associative, commutative, flexible (that is, $x (yx) = (xy)x$, for
	$x, y\in E$); not necessarily power-associative (that is, the subalgebra generated by a single element can not be associative); and
	preserved by direct sums. Moreover, evolution algebras are not closed under subalgebras. So, a subalgebra of $E$ is called an {\em evolution subalgebra} if the former is spanned by a subset of generators of the latter. As such, it is an ideal of $E$. The latter is {\em indecomposable} if it is not the direct sum of two nonzero ideals; {\em connected} if it is not the direct sum of two proper evolution subalgebras; {\em irreducible} if it has no proper subalgebra; and {\em simple} if it has no proper evolution subalgebra. Tian proved \cite[Theorem 9]{Tian2008} that any finite-dimensional evolution algebra has a simple evolution subalgebra.
	
	Every evolution algebra $E$ is uniquely determined by its {\em evolution operator} $L:E\mapsto E$, which is linearly described from
	\[L(e_i):=e_ie_i =  \sum_i p_{ik} e_k, \text{ for all } i\in\{1,\ldots,n\}.\]
	Then, for each positive integer $m>2$, it is defined the {\em plenary power}
	\[e_i^{[m]}:=L\left(e_i^{[m-1]}\right),\]
	where $e_i^{[0]}:=e_i$. The operator $L$ may be considered as time-step in a discrete-time dynamical system, and hence, it describes the dynamical flow of the evolutionary process represented by the evolution algebra $E$.
	In this regard, a generator $e_i$ of the evolution algebra $E$ is called {\em algebraically persistent} if the evolution subalgebra generated by $e_i$ is a simple evolution subalgebra. As such, it represents an absorbent state giving rise to a stable evolution. Otherwise, it is said to be {\em algebraically transient}, which gives rise to a transitory evolution. The latter would eventually derive to an absorbent state, but only after certain generations. Any generator of a simple evolution subalgebra is algebraically persistent. In particular, every finite-dimensional evolution algebra contains an algebraically persistent generator. 
	
	A connected evolution algebra is simple if and only if all its generators are algebraically persistent. Tian proved \cite[Theorem 11]{Tian2008} that, if $E$ is a connected finite-dimensional evolution algebra, then
	
	\[E=\bigoplus_{j=1}^{n_0} E_{0,j} + B_0,\]
	where $\bigoplus$ and $+$ denote, respectively, direct sums of either subalgebras or linear subspaces; each $E_{0,j}$ is a simple evolution subalgebra so that $E_{0,j}\cap E_{0,j'}=\{0\}$, whenever $j\neq j'$; and $B_0$ is a linear subspace spanned by algebraically transient generators, which is called {\em transient space}. Even if $B_0$ is not an evolution subalgebra in general, it may be endowed of evolution structure from $E$ once the multiplication is restricted within $B_0$. The previous decomposition process can inductively be repeated until no transient space arises. That is, there exist an integer $m>0$ and a subset $\{n_1,\ldots,n_m\}\subset\mathbb{N}$ such that
	\[E=\sum_{i=0}^m\bigoplus_{j=1}^{n_i} E_{i,j},\]
	where $\sum$ denote direct sum of linear subspaces and each $E_{i,j}$ is a simple evolution algebra so that $E_{i,j}\cap E_{i,j'}=\{0\}$, whenever $j\neq j'$. Each subset $\{E_{i,1},\ldots,E_{i,n_i}\}$ constitutes the $i^{\mathrm{th}}$ level of a hierarchical structure describing the dynamical flow of $E$, which enabled Tian \cite{invita} to ensure that  {\em ``simple evolution algebras are the basic blocks for building general evolution algebras''}. Tian dealt in particular with the spectrum of the evolution operator $L$ at the $0^{\mathrm{th}}$ level of this hierarchy, and set out the following problem.
	
	\begin{prob}[\cite{Tian2008}]\label{prob_spectrum} Study the spectrum of the evolution operator of an evolution algebra at the $i^{\mathrm{th}}$ level of the hierarchical structure of the latter, for each $i>1$. Furthermore, study the spectrum theory concerning plenary powers of evolution algebras.
	\end{prob}
	
	Furthermore, the just described hierarchical structure may also be used to classify evolution algebras. More specifically, evolution algebras can be distributed according to the number of levels of its related hierarchical structure and the number of simple evolution subalgebras at each level. Every evolution algebra is homomorphic to a unique evolution algebra, which is called its {\em skeleton-shape}, whose hierarchy has only one-dimensional evolution algebras in all its levels \cite[Theorem 15]{Tian2008}. In particular, evolution algebras within each class of the mentioned distribution are uniquely determined up to their skeleton-shape homomorphism.
	
	Another method for classifying evolution algebras is based on the period of their generators. In this regard, a generator $e_i$ is said to {\em occur} in $x=\sum_{i=1} ^n \alpha_ie_i\in E$ if $\alpha_i\neq 0$. It is denoted $e_i\prec x$. The {\em period} of a generator $e_i$ is defined as
	\[d_i:=g.c.d.\left\{\mathrm{log}_2 m\colon\, e_i\prec e_i^{[m]}\right\}.\]
	If $\left\{\mathrm{log}_2 m\colon\, e_i\prec e_i^{[m]}\right\}=\emptyset$, then $d_i=\infty$. Further, if $d_i=1$, then the ge\-ne\-ra\-tor $e_i$ is called {\em aperiodic}. Otherwise, it is called {\em periodic}. If the evolution algebra is non-negative simple, then all its generators have the same period \cite[Theorem 7]{Tian2008}. As such, non-negative simple evolution algebras may be classified as either periodic or aperiodic. In the first case, if $d>1$ is the period of all the generators of a non-negative simple evolution algebra $E$, then these generators can be partitioned into $d$ disjoint classes $C_0,\ldots,C_{d-1}$ such that, if $\Delta_i$ denotes the linear subspace spanned by $C_i$, for all $i\in \{0,\ldots,d-1\}$, then $L(\Delta_i)\subseteq \Delta_{i+1 \mod d}$, and
	
	\[E=\bigoplus_{i=0}^{d-1}\Delta_i.\]
	
	Based on the already mentioned fact that simple evolution algebras may be seen as the basic blocks on the hierarchical structure of any evolution algebra, Tian proposed to generalize the previous results on non-negative simple evolution algebras. In particular, he set out the following problem.
	
	\begin{prob}[\cite{invita}]\label{prob_simple} Delve into the general study of simple evolution algebras, and not only non-negative ones. Characterize them and prove the existence of transitive occurrence relations. To this end, study in particular simple evolution algebras over finite fields. Furthermore, it is also required to delve into the study of simple evolution algebras at higher levels of the hierarchy of evolution algebras, and not only at the $0^{\mathrm{th}}$ level.
	\end{prob}
	
	\section{Developing the theory of evolution algebras}\label{sec:FirstSteps}
	
	This section deals with some of the most important results concerning the mathematical fundamentals of the theory of evolution algebras. Again, let $E$ denote from now on an $n$-dimensional evolution algebra over a base field $\mathbb{K}$, with natural basis $\{e_1,\ldots,e_n\}$ and structural constants $p_{ik}\in\mathbb{K}$, for all $i,k\in\{1,\ldots,n\}$.
	
	\subsection{Nilpotent evolution algebras}
	
	In 2013, based on the non-necessary associativity of evolution algebras, Jos\'e Manuel Casas, Ladra, Omirov and Rozikov \cite{Casas2013} distinguished among nil evolution algebras, right nilpotent evolution algebras and nilpotent evolution algebras. More specifically, for each vector $x\in E$ and each positive integer $i>1$, let $x^{<i>}:=x^{<i-1>}x$, where $x^{<1>}:=x$. Then, the vector $x$ is {\em nil} if there exists a positive integer $n$ such that $x^{<n>}=0$. The evolution algebra $E$ is {\em nil} if all its elements are nil. Its {\em nilindex} is the minimum positive integer $n$ such that $x^{<n>}=0$, for all $x\in E$. Further, the evolution algebra $E$ is {\em right nilpotent} if $E^{<n>} = 0$, for some positive integer $n$, where $E^{<1>}:=E$ and $E^{<k+1>}:= E^{<k>} \, E$, for every positive integer $k>1$.  Finally, the evolution algebra $E$ is {\em nilpotent} if $E^n = 0$, for some positive integer $n$, where $E^1:=E$ and $E^k:=\sum_{i=1}^{k-1} E^i E^{k-i}$, for every positive integer $k>1$. The index of (right) nilpotency is defined similarly to the nilindex.
	
	This distinction is not necessary in case of dealing with finite-dimensional complex evolution algebras. The mentioned authors proved \cite[Theorem 3.1]{Casas2013} the equivalence of right nilpotency and nility for finite-dimensional evolution algebras. Their equivalence with nilponency was shortly after proved by Luisa Mar\'\i a  Camacho, Jos\'e Ram\'on G\'omez, Omirov and Rustam Turdibaev \cite[Corollary 4.6]{LGOT} in case of dealing with finite-dimensional complex evolution algebras. Their structure matrices are all of them upper triangular after a suitable permutation of their corresponding natural bases (see also \cite[Theorem 2.7]{CLOR}). In any case, it was observed that the nilindex and the indexes of (right) nilpotency of one such an evolution algebra do not coincide in general. Further, it was proved \cite[Theorem 4.5]{LGOT} that the maximum index of nilpotency of an $n$-dimensional evolution algebra is $2^{n-1}+1$. A characterization of evolution algebras reaching this maximum value was described in \cite[Theorem 3.6]{Casas2013}.  More specifically, the index of nilpotency of our evolution algebra $E$ is $2^{n-1}+1$ if and only if $p_{12}p_{23}\ldots p_{(n-1)n}\neq 0$.  The automorphism group and local derivation of these evolution algebras with maximum index of nilpotency was studied by Mukhamedov, Otabek Khakimov, Omirov and Qaralleh \cite{Mukhamedov2019}. Their distribution into isomorphism classes over the complex field was established in \cite{Mukhamedov2020}. Evolution algebras of index of nilpotency $2^{n-2}+1$ were similarly studied by Qaralleh \cite{Qaralleh2020}.
	
	Also in 2013, Ahmed Sadeq Hegazi and Hani Abdelwahab \cite{Hegazi2013} classified nilpotent evolution algebras of dimension up to three over any base field, and of dimension four for algebraically closed fields of any characteristic and also over the real field. To this end, they studied the annihilator extensions of evolution algebras. Recall that the {\em annihilator} of an algebra $A$ is the set
	
	\[\mathrm{Ann}
	(A)=\{x\in A\colon\, xy=0, \text{ for all } y\in A\}.\]
	According to Alberto Elduque and Alicia Labra \cite[Lemma 2.7]{Ella2}, the annihilator of the evolution algebra $E$ is $\mathrm{Ann}
	(E)=\langle\,e_i\in A\colon\, e_i^2=0\,\rangle$.
	
	In 2014, Tian and Yi Zou \cite{TZ} characterized nil evolution algebras and nilpotent evolution algebras in terms of their structural constants. In addition, they constructed finitely generated nil evolution algebras that are not nilpotent and explained how their results may be used to model certain population dynamics. They also studied  finite-dimensional complex evolution algebras \cite{LGOT}, for which they characterized those nilpotent evolution algebras that are isomorphic to evolution algebras in Jordan normal form. In 2015, Abror Khudoyberdiyev, Omirov and Qaralleh \cite{KOQ2015} reduced this study to idempotent and absolute nilpotent elements, whose dynamic in chains of two-dimensional evolution algebras depending on the time  was dealt with by Sherzod Murodov \cite{Muro2}.
	
	In 2016, Elduque and Labra \cite{Ella} distributed into isomorphism classes the four- and five-dimensional indecomposable nilpotent evolution algebras over algebraically closed fields of characteristic different from two. One year later, Abdelwahab Alsarayreh, Qaralleh and Muhammad Ahmad characterized \cite{Alsarayreh2017} the space of derivations of three-dimensional nilpotent evolution algebras. In 2018, Camacho, Khudoyberdiyev and Omirov \cite{CKO} distributed into isomorphism classes those nilpotent evolution algebras for which any subalgebra is an evolution subalgebra. Shortly after, Omirov, Rozikov and Mar\'\i a Victoria Velasco \cite{ORVlast} distributed into isomorphism classes the set of nilpotent evolution algebras of dimension $n\leq 5$.
	
	\subsection{Power-associative evolution algebras}
	
	Evolution algebras are not necessarily power-associative. Yolanda Ca\-brera, Mercedes Siles and Velasco \cite{CSV16} observed that the evolution algebra $E$ is power-associative only if $p_{ii}^2=p_{ii}$, for all $i$. Thus, every nil evolution algebra is power-associative \cite[Proof of Theorem 2.2]{CLOR}. Examples of nil associative evolution algebras, which are therefore power-associative, are shown in \cite{Ella}. Further, Luisa Mar\'\i a  Camacho, Jos\'e Ram\'on G\'omez, Omirov and Rustam Turdibaev described \cite[Example 4.8]{LGOT} a four-dimensional evolution algebra that is power-associative, but it is not associative.
	
	\newpage
	
	In 2020, Moussa Ouattara and Souleymane Savadogo proved \cite[Theorem 6]{Outtara2020a} that an evolution algebra is power-associative if and only if it is a Jordan algebra. They also determined all the power-associative evolution algebras up to dimension six. The same authors also dealt \cite{Outtara2020} with the distribution into isomorphism classes of those indecomposable and not power-associative evolution nil-algebras up to dimension five, having nil-index four.
	
	\subsection{Baric evolution algebras}
	
	Unlike genetic algebras, evolution algebras are not baric in general. In this regard, Casas, Ladra and Rozikov proved  \cite[Theorem 3.2]{Casas2011} that an $n$-dimensional real evolution algebra is baric if and only if there exists a positive integer $k\leq n$ such that $e_ke_k=p_{kk}e_k\neq 0$. Its weight $\omega$ is defined so that $\omega\left(\sum_i a_i e_i\right)=p_{kk}a_i$. Much more recently, Ouattara and Savadogo \cite[Theorem 3.3]{Outtara2020} have shown that this result can be extended to any commutative field of characteristic distinct from two. Then, they have also proved \cite[Corollary 3.4]{Outtara2020} that every baric evolution algebra $E$ over a base field $\mathbb{K}$ admits a natural basis $\{e_1,\ldots,e_n\}$ and a weight $\omega:E\rightarrow \mathbb{K}$ such that $\omega(e_1)=1$ and $\omega(e_i)=0$, for all $i>1$. Moreover, $E=\mathbb{K}e_1\oplus \mathrm{ker}(\omega)$, with $e_1\mathrm{ker}(\omega)=0$. 
	
	In this last work, it was also proved the non-existence of evolution train algebras of rank two \cite[Proposition 3.7]{Outtara2020} and the fact that a baric evolution algebra of weight $\omega$ is a train algebra of rank $r+1>2$ if and only if $\mathrm{ker}(\omega)$ is nil of nil-index $r$ \cite[Theorem 3.8]{Outtara2020}. Moreover, any evolution train algebra is a special train algebra \cite[Theorem 3.12]{Outtara2020}.
	
	Casas, Ladra, Omirov and Rozikov proved \cite[Theorem 4.1]{Casas2013} that finite-dimensional nilpotent evolution algebras are not dibaric. They also proved \cite[Corollary 4.5]{Casas2013} that an evolution algebra is not dibaric when the determinant of its structure matrix is not zero. Finally, they characterized \cite[Proposition 4.10]{Casas2013} the two-dimensional real dibaric evolution algebras.
	
	In 2021, Andr\'e Conseibo, Savadogo and Ouattara \cite{Conseibo2021} have characterized baric evolution algebras that are Bernstein algebras and have classified into isomorphism classes these type of algebras, for dimension $n\leq 4$.
	
	\subsection{(Semi)simple evolution algebras}
	
	In 2016, Cabrera, Siles and Velasco established \cite[Corollary 4.6]{CSV16} that a finite-dimensional evolution algebra is simple if and only if its structure matrix associated to a given natural basis is nonsingular and that basis cannot be reordered so that, for some $m<n$, the structure matrix has block form
	
	\[\left(\begin{array}{cc}
		W_{m\times m} & U_{m\times (n-m)}\\
		0_{(n-m)\times m} & Y_{(n-m)\times (n-m)}
	\end{array}\right).\]
	
	In 2019, Velasco \cite[Corollary 3.7]{Victoria} characterized the maximal modular ideals of an evolution algebra and proved \cite[Corollary 3.8]{Victoria} that all of them have codimension one. Recall in this regard that a {\em modular ideal} of the evolution algebra $E$ is any ideal $M$ of the algebra that is endowed with a {\em modular unit} $u\in E$ such that $x-xu\in M$, for all $x\in E$. In particular, Velasco \cite[Proposition 2.2]{Victoria} proved that a finite-dimensional evolution algebra has a unit (and hence, all its ideals are modular) if and only if its structure matrix associated to a given natural basis is diagonal with non-zero entries.
	
	The intersection of all the maximal modular ideals of the algebra constitutes its {\em Jacobson radical}. Then, an evolution algebra is called {\em semisimple} if its Jacobson radical is zero. Velasco dealt with semisimple evolution algebras by introducing the so-called {\em spectrum} and {\em $m$-spectrum} of an evolution algebra. Both concepts were described in terms of the eigenvalues of a suitable matrix related to the structure matrix of the algebra. She set that the algebra is {\em $m$-semisimple} (respectively, {\em spectrally semisimple}) if the zero ideal is the only one for which the $m$-spectrum is $\{0\}$ (respectively, the spectrum is $\{0\}$). Unlike the associative case (where semisimplicity, spectrally semisimplicity and $m$-semisimplicity are equivalent), there exist $m$-semisimple evolution algebras whose Jacobson radical coincides with the algebra.
	
	\subsection{Perfect evolution algebras}
	
	An algebra $A$ is {\em perfect} if $A^2=A$. From \cite[Theorem 4.4]{Ella2}, every perfect evolution algebra has a unique natural basis. (Nadia Boudi, Cabrera and Siles found  \cite[Corollary 2.7]{Boudi2020} a condition for an evolution algebra to have a unique natural basis.) Moreover, the automorphism group of a finite-dimensional perfect evolution algebra is finite \cite[Theorem 4.8]{Ella2}.
	
	In 2019, Cabrera, M\"uge Kanuni and Siles \cite{CKS} introduced the {\em basic ideal} of an evolution algebra as an ideal having a natural basis that can be extended to a natural basis of the whole algebra. They established that maximal basic ideals are unique except for those ideals having codimension one. This fact enabled them to ensure \cite[Proposition 2.19]{CKS} that the number of zeros in the structure matrix of an evolution algebra is an invariant, except for those cases in which there exists a maximal basic ideal of codimension one.
	
	It is so that the authors showed how maximal basic ideals play a fundamental role in the distribution into isomorphism classes of four-dimensional perfect non-simple evolution algebras. They dealt with this classification over a field with mild restrictions. The complete distribution into isomorphism classes was later studied in \cite{Behn2020}.
	
	That one of two dimensional perfect evolution algebras over arbitrary fields has recently been established in \cite{Cardoso2021}, where it was also studied the automorphism group, the space of derivations and identities of degree at most four of both perfect and non-perfect evolution algebras.
	
	In 2020, Boudi, Cabrera and Siles proved \cite[Proposition 4.2]{Boudi2020} that every nonzero ideal in a perfect evolution algebra is basic. In addition, they observed \cite[Corollary 3.8]{Boudi2020} that a perfect evolution algebra has a nilpotent element of order three only if its structure matrix associated to a given natural basis has a vanishing principal minor. The reciprocal also holds when every element of the base field is a square.
	
	In 2021, Elduque and Labra \cite[Theorem 4.1]{Elduque2021} have characterized the space of derivations of any perfect evolution algebra. Particularly, they have proved that, if the characteristic of the base field is zero or two, then this subspace is only formed by the null map. Otherwise, its dimension depends on the number of connected components of a related graph.
	
	\subsection{Derivation of evolution algebras}
	
	In the previous subsections, we have already mentioned some works concerning the study of the space of derivations of different types of evolution algebras. Let us recall here that a {\em derivation} of the evolution algebra $E$ is any linear map $d:E\rightarrow E$ such that $d(xy)=d(x)y+xd(y)$, for all $x,y\in E$. 
	
	As for any algebra, the space of derivations of a given algebra is a Lie algebra that constitutes a good approximation to its automorphism group and hence, to its algebraic structure. The system of equations described by this space was given by Tian himself \cite{Tian2008}. More specifically, if $d$ is a derivation of the evolution algebra $E$ that is described so that $d(e_i)=\sum_k d_{ik}e_k$, with $d_{ik}\in \mathbb{K}$, for all $i$, then it must be 
	\[p_{ik}d_{ji}+p_{jk}d_{ij}=0, \text{ for all } i,j,k \text{ such that } i\neq j,\]
	and
	\[\sum_k p_{ik}d_{kj}=2p_{ij}d_{ii},\text{ for all } i,j.\]
	
	In 2013, Camacho, G\'omez, Omirov and Turdibaev \cite{Camacho2013} studied formally the space of derivations of $n$-dimensional complex evolution algebras depending on the rank of certain matrices. In particular, they proved that such a subspace is zero for evolution algebras having a non-singular structure matrix \cite[Theorem 2.1]{Camacho2013}. They also described the space of derivations of those $n$-dimensional evolution algebras having a structure matrix of rank $n-1$.
	
	Much more recently, Paula Cadavid, Mary Luz Rodi\~no and Pablo M. Rodr\'\i guez \cite{cadavid2} described the space of derivations of those evolution algebras that are uniquely associated with finite simple and connected graphs of order $n\geq 3$ according to Tian's proposal (see Subsection \ref{subsec:Graph}). It was observed that each one of these spaces of derivations depends on the twin partition of the set of vertices of the corresponding graph. (Recall here that two vertices within a graph are called {\em twins} if their neighborhoods coincide.) In particular, if every part in this twin partition contains at most two vertices, then the space of derivations is only formed by the null map \cite[Theorem 2.3]{cadavid2}. The space of derivations in case of existing some part within this twin partition with at least three vertices was also characterized \cite[Theorem 2.6]{cadavid2}. It is remarkable the fact that the authors dealt with examples of finite-dimensional evolution algebras having structure matrices of any rank. 
	
	The same authors, together with Cabrera, have recently explored \cite{Cabrera2021} a similar extended approach by considering to this end the non-zero entries of the structure matrix of the evolution algebra under consideration and its associated directed graph, which had previously been described by Cabrera, Siles and Velasco \cite{CSV16}. By making use of this approach, the authors have characterized the derivations of non-degenerate irreducible three-dimensional evolution algebras.
	
	Let us finish this subsection with a proposal of further work on the topic under consideration. Similarly to the relationship between the subset of derivations of an algebra and its automorphism group, there also exists a relationship between its {\em autotopism group} (that is, the group of isotopisms preserving the algebra $A$ under consideration) and its subset of ternary derivations. This terminology was introduced for general algebras by Clara Jim\'enez-Gestal and Jos\'e Mar\'\i a P\'erez-Izquierdo \cite{Jimenez2003}. More specifically, a {\em ternary derivation} of an algebra $A$ is any triple $(d_1,d_2,d_3)$ of endomorphisms of the algebra such that 
	\[d_1(xy)=d_2(x)y + xd_3(y),\]
	for all $x,y\in A$. (See \cite{Barquero2021} for a comprehensive motivation of this concept.) The following question arises from the important role that isotopisms play as algebraic representations of mutations in Genetics.
	
	\begin{prob} Characterize the subspace of ternary derivations of any given finite-dimensional evolution algebra and establish the relationship with its autotopism group.
	\end{prob}
	
	\subsection{Classification of evolution algebras}
	
	In the previous subsections, the distribution into isomorphism classes of different types of evolution algebras has been outlined. Let us indicate here some more results concerning the classification of finite-dimensional evolution algebras.
	
	In 2013, Casas, Ladra, Omirov and Rozikov  distributed \cite[Theorem 4.1]{CLOR} all the two-dimensional  evolution algebras over the complex field into six non-isomorphic classes. Nevertheless, Cabrera, Siles and Velasco \cite{CSV} realized that this classification did not consider the evolution algebra with natural basis $\{e_1,e_2\}$ verifying $e_1^2=e_2$ and $e_2^2=e_1$. It was noticed in the classification of  three-dimensional evolution algebras with two-dimensional derived algebra and one-dimensional annihilator \cite[Theorem 3.5]{CSV}.
	
	In 2014, Murodov distributed \cite[Theorem 1]{Murodov2014}  into isomorphism classes the set of two-dimensional evolution algebras over the real field. More recently, the classification over a general field was independently achieved by \'Oscar J. Falc\'on, Ra\'ul M. Falc\'on and Juan N\'u\~nez \cite{FFN2}, and by Maria Inez Cardoso, Daniel Gon\c{c}alves, Dolores Mart\'\i n Barquero, C\'andido Mart\'\i n Gonz\'alez and Siles \cite{Cardoso2021}. In \cite{FFN2}, it was also observed that the distribution of evolution algebras into isotopism classes is uniquely related with the mutation of alleles in non-Mendelian Genetics. In that reference, it was obtained the distribution into four isotopism classes of all two-dimensional evolution algebras, whatever the base field is.
	
	In 2017, Cabrera, Siles and Velasco \cite{CSV} classified into  116 distinct types the set of three-dimensional evolution algebras over any field of characteristic distinct from two and in which there are roots of orders two, three and seven. Their results agreed with those ones of Elduque and Labra \cite{Ella} concerning the classification of indecomposable nilpotent evolution algebras of dimension up to five over any algebraically closed fields of characteristic distinct from two. More recently, based on an alternative approach, Mar\'\i a Eugenia Celorrio and Velasco \cite[Theorem 11]{Celorrio2019} classified the mentioned 116 types of algebras into 14 non-isomorphic types,
	
	In 2018, the distribution of three-dimensional evolution algebras having one-dimensional annihilator was determined in \cite[Theorem 4.7]{FFN}, where it was also established \cite[Corollary 3.3]{FFN} the distribution of three-dimensional evolution algebras into isotopism classes, whatever the base field is. This last classification enables one to describe the spectrum of genetic patterns of three distinct genotypes during a mitosis process. 
	
	Also in 2018, Anvar N. Imomkulov \cite{Imon} introduced an evolution operator for evolution algebras. Both the set of fixed points and the Jacobian matrix of this operator were studied for two-dimensional evolution algebras. In particular, it was established the distribution into isomorphism classes of those two-dimensional evolution algebras having this Jacobian matrix as structure matrix. The distribution of the three-dimensional case was also dealt with by the same author \cite{Imon2019}, who proposed as further work the following research problem.
	
	\begin{prob}\label{problem_Approx} Study possible approximations of finite-dimensional algebras by means of evolution algebras.
	\end{prob}
	
	A first approach of this problem was considered by Imomkulov himself, together with Rozikov, who determined \cite{Imon2019a} in 2020 the relationship among two- and three-dimensional Leibniz algebras and evolution algebras. Also in 2020, Imomkulov \cite{Imomkulov2020} got a result about absolute nilpotent elements of evolution algebras corresponding to approximation of finite dimensional algebras. Already in 2021, Cabrera et al. \cite{Cabrera2021a} have provided four different constructions producing three-dimensional evolution algebras from two-dimensional algebras.
	
	\section{Relationship to other topics}\label{sec:Topics}
	
	This section outlines the relationship among evolution algebras, Graph theory, Group theory, Markov chains and Biology.
	
	\subsection{Graph theory}\label{subsec:Graph}
	
	In his original manuscript, Tian \cite{Tian2008} realized that every finite graph $G=(V,E)$ may be associated with an evolution algebra $A(G)$ having the vertex set $V=\{e_1,\ldots,e_n\}$ as its natural basis, and such that
	\[e_ie_i=\sum_{e_k\in\Gamma(e_i)}e_k,\]
	for every positive integer $i\leq n$, where $\Gamma(e_i)$ is the neighbourhood of the vertex $e_i$ in the graph $G$. In this way, isomorphic graphs give rise to isomorphic evolution algebras \cite[Theorem 41]{Tian2008}. Moreover, if
	\[L^m(e_i)=\sum_{j=1}^n p_{ik}e_k,\]
	where $L$ is the evolution operator associated to $A(G)$, then $p_{ik}$ coincides with the total number of paths of length $m$ between both vertices $e_i$ and $e_j$ in the graph $G$ \cite[Theorem 43]{Tian2008}.
	
	In order to deal with the reciprocal, Tian \cite{Tian2008} introduced the concept of {\em graphicable algebra} as an evolution algebra of natural basis $V=\{e_1, \dots, e_n\}$ such that
	\[e_ie_i = \sum_{e_k\in V_i} e_k,\]
	for all positive integer $i\leq n$,  where $V_i\subseteq V$. Even if evolution algebras are not graphicable in general, every graphicable
	algebra is uniquely associated to a graph such that $V_i=\Gamma(e_i)$. In this way, two isomorphic graphicable
	algebras give rise to isomorphic graphs \cite[Theorem 42]{Tian2008}. In particular, Tian introduced {\em cycle algebras}, {\em path algebras} and {\em complete algebras} as graphicable algebras giving rise, respectively, to cycles, paths and complete graphs. In a similar way, N\'u\~nez, Mar\'\i a Luisa Rodr\'\i guez-Ar\'evalo and Mar\'\i a Trinidad Villar \cite{NRV} dealt with complete tripartite and $n$-partite, star, friendship, wheel and snark graphicable algebras. These authors also proved \cite[Theorem 3.8]{NRV} the existence of graphicable subalgebras whithin a graphicable algebra. Further, N\'u\~nez, Marithania Silvero and Villar dealt \cite{AMC} with the particular case of graphicable algebras for which $e_k\in V_i$ if and only if $e_i\in V_k$, and $e_i\not\in V_i$, for every pair of positive integers $i,k\leq n$ such that $i\neq k$. They called {$S$-graphicable algebras} this type of algebras.
	
	Based on both notions of evolution algebra arising from a graph and graphicable algebra, Tian remarked \cite{Tian2008} that the {\em `intrinsic and coherent relation of evolution algebras with graph theory allows to analyze graphs algebraically [...] and graph theory may be used as a tool to study non-associative algebras''}. It is so that he asked whether {\em ``every statement or problem in graph theory can be translated into the language of evolution algebras''}. He put particular interest in the possible relationship of evolution algebras with random graphs and networks, random walks on graphs, weighted graphs and directed graphs.
	
	Furthermore, Tian also asked for those evolution algebras arising from finite graphs whose associated Ihara-Selberg zeta function satisfies the Sunada's analogue of Riemann hypothesis \cite{Sunada1986}.
	
	In 2011, Rozikov and Tian \cite{Rozikov2011} described an alternative for relating an evolution algebra to a given finite and simple graph by means of state spaces endowed with Gibbs measures. For connected graphs, they determined the hierarchical structure of these algebras, together with the dimension and number of one- and four-dimensional evolution subalgebras \cite[Theorem 3.2]{Rozikov2011}. Moreover, all the evolution algebras arising from the same finite and simple graph, but defined by different Gibbs measures, are pairwise isomorphic \cite[Theorem 3.3]{Rozikov2011}. The authors established as further open problem the study of the hierarchical structure of evolution algebras associated to graphs that are either finite, but not connected; or countable and connected.
	
	In 2015, Elduque and Labra \cite{Ella2} defined the {\em directed graph attached} to an $n$-dimensional evolution algebra of structure matrix $(p_{ik})_{i,k}$ with respect to a given natural basis as the graph of set of vertices $V=\{1,\ldots,n\}$ and set of arcs $A=\{(i,k)\in V\times V\colon\, p_{ik}\neq 0\}$. If each arc $(i,k)\in A$ is labeled as $p_{ik}$, then one obtains the {\em directed weighted graph attached} to $E$. If $\mathrm{Ann}(E)=\{0\}$, then its directed graph is connected if and only if $E$ does not split into a direct sum of simple ideals \cite[Proposition 2.8]{Ella2}. If the annihilator is not trivial, then this decomposition is not possible if and only if its directed graph is connected for every natural basis \cite[Proposition 2.10]{Ella2} (see also \cite[Proposition 5.4]{CSV16}). Further, $E$ is nil if and only if its attached directed graph contains no oriented cycles \cite[Theorem 3.4]{Ella2}. These attached (weighted) directed graphs were also used to distribute into isomorphism classes the set of four- and five-dimensional nilpotent evolution algebras \cite{Ella}.
	
	In 2020, Cadavid, Rodi\~no and Rodr\'\i guez \cite{cadavid} described the evolution algebra induced by the random walk on a graph and studied its relationship with the evolution algebra determined by the same graph according to Tian's proposal \cite{Tian2008}. Recall here that the {\em random walk} on a graph of set of vertices $V$ and adjacency matrix $(a_{ij})_{i,j}$ is a discrete Markov chain with state space $V$ such that the transition probability of moving from  state $i$ to state $j$ is
	\[p_{ij}=\frac {a_{ij}}{\sum_{k\in V} a_{ik}}.\]
	The mentioned authors described the {\em evolution algebra induced} by this random walk as the evolution algebra of natural basis $\{e_1,\ldots,e_n\}$ and structure matrix $(p_{ij})_{i,j}$.
	
	Then, they studied under which conditions one may ensure that this evolution algebra is isomorphic to Tian's evolution algebra associated to the graph under consideration. More recently, the same authors have proved \cite[Theorem 2.2]{Cadavid2021} that both algebras are indeed strongly isotopic. In addition, they have provided conditions \cite[Theorem 2.3]{Cadavid2021} under which these algebras are isomorphic or not.
	
	Also in 2020, Manuel Ceballos, N\'u\~nez and \'Angel Tenorio \cite{CNT} described the {\em directed weighted pseudograph associated} to an $n$-dimensional evolution algebra of structure matrix $(p_{ik})_{i,k}$ as the pseudograph resulting after adding a loop on each vertex $i$ of the Elduque and Labra's directed weighted graph attached to $E$, whenever $p_{ii}\neq 0$. This loop is then labeled as $p_{ii}$. The distribution into isomorphism classes of these pseudodigraphs enabled them to establish a new classification of evolution algebras.
	
	Furthermore, Rafael Gonz\'alez and N\'u\~nez \cite{Gonzalez2020} translated into the algebraic language some basic concepts and results concerning the directed graphs associated to an evolution algebra. In particular, based on the adjacency of graphs, they introduced the notions of adjacency, walk, trail, circuit, path and cycle of an evolution algebra. In addition, strongly and weakly connected evolution algebras were introduced as the algebraic equivalences of the same concepts in graph theory. It enabled the authors to introduce the notions of distance, girth, circumference, eccentricity, center, radio, diameter and geodesic of an evolution algebra, together with the concepts of Eulerian and Hamiltonian evolution algebras.
	
	In 2021, Qaralleh and Mukhamedov \cite{Qaralleh2021} have introduced the concept of {\em Volterra evolution algebra} as an evolution algebra whose structural matrix is described by skew symmetric matrices. This type of algebras are not nilpotent. In particular, they have proved \cite[Theorem 6.5]{Qaralleh2021} that the directed weighted graphs associated to two Volterra evolution algebras are isomorphic if and only if the mentioned algebras are isomorphic.
	
	\subsection{Group theory} Tian \cite{Tian2008} associated any given group $(G,\circ)$ having a finite set of generators $S$ with the evolution algebra described over a base field $\mathbb{K}$ so that, for each $g\in G$, one has
	\begin{equation}\label{eq:EA_group}
		g\cdot g =\sum_{e\in S}k_eg\circ e,
	\end{equation}
	for some $k_e\in\mathbb{K}$. Particularly, Tian set out the following question.
	
	\begin{prob}[\cite{Tian2008}]\label{prob_EA_groups} How can the properties of a group be translated to the corresponding evolution algebra?
	\end{prob}
	
	\subsection{Markov chains}
	
	Tian himself \cite{Tian2008} realized that every Markov chain is related to a Markov evolution algebra whose structural constants coincide with the transition probabilities of the former and whose basis constitutes the state space of the Markov process.
	
	As such, generators of Markov evolution algebras represent states of stochastic processes described by Markov chains, which may, therefore, be studied by means of the theory of evolution algebras. In particular, Markov chains may be classified according to the hierarchies of their corresponding evolution algebras. Tian proved that a Markov chain is irreducible if and only if its related evolution algebra is simple \cite[Theorem 18]{Tian2008}, and that a subset of state space of a Markov chain is probabilistic closed if and only if it generates an evolution subalgebra \cite[Theorem 17]{Tian2008}. In addition, probabilistic periods in Markov chains are uniquely related to periods of generators in evolution algebras \cite[Proposition 10]{Tian2008}.
	
	Tian also introduced the possibility of dealing with {\em continuous evolution algebras}, once structural constants are replaced by differential functions depending on a variable. If the latter represents the time, then continuous evolution algebras are uniquely related to continuous-time Markov chains \cite{TianXiao2006}. Tian thought that the study of continuous evolution algebras {\em ``will be very interesting because they have a kind of semi-Lie group structure''}.
	
	In 2011, as a generalization of Tian's idea, Casas and Rozikov \cite{Casas2011} realized that every Markov chain is associated to a {\em chain of $n$-dimensional evolution algebras} $\{A^{[s,t]}\colon\,0\leq s\leq t\}$ depending on two time parameters $s$ and $t$. Their structure matrices $\mathcal{M}^{[s,t]}$ are all of them stochastic and satisfy that 
	\[\mathcal{M}^{[s,t]}= \mathcal{M}^{[s,\tau]}\mathcal{M}^{[\tau,t]}\] for all $s<\tau<t$. This chain represents a continuous time dynamical system that is an evolution algebra in each fixed time. If all the matrices of the chain coincide, then they constitute the same evolution algebra that Tian associated to a Markov chain. Furthermore, if both time parameters $s$ and $t$ are reduced to only one of them, then the chain is called {\em time-homogeneous} and is associated to a time-homogeneous Markov chain. The authors also dealt with chains of evolution algebras for which the stochastic condition is removed.
	
	In 2013, Rozikov and Mudorov  \cite{RM2} constructed $25$ distinct examples of chains of two-dimensional evolution algebras, for which they studied their baricity, nilpotency and idempotency. (Some new examples have been recently described by Ladra and Murodov \cite{Ladra2021}.) In 2020, Imomkulov and Velasco \cite{ImomVelasco2020} got a description depending on rank of chains of three-dimensional evolution algebras. And then, for these chains, Imomkulov \cite{Imomkulov2020a} studied the behavior of both the baric property and the set of absolute nilpotent elements, and also the time-depending dynamics of the set of idempotent elements.
	
	The general case was already considered in 2015, when Omirov, Rozikov and Kaisar Tulenbayev \cite{ORT2015} dealt with real chains of $n$-dimensional evolution algebras corresponding to a permutation of $n$ numbers and showed that one such a chain is trivial if and only if the permutation has no fixed points \cite[Proposition 1]{ORT2015}. (The study of evolution algebras arising from permutations was introduced by Narkuziev \cite{Narkuziev2014}, who described a nilpotency condition on this type of algebras.) Moreover, since every trivial chain of evolution algebras is a chain of nilpotent evolution algebras \cite[Proposition 2]{ORT2015}, they constructed chains of three-dimensional evolution algebras. They also dealt with the construction of arbitrary-dimensional symmetric chains of evolution algebras.
	
	In 2017, Ladra and Rozikov \cite{Ladra2017, Ladra2017a} described non-homogeneous continuous time Markov chains. More recently, Irene Paniello \cite{Paniello2020} has delved into the algebraic structure of Markov evolution algebras for both the discrete-time and the continuous-time arising from standard stochastic semigroups.
	
	\subsection{Biology}
	
	Since the original manuscripts of Tian and Vojtechovsky, evolution algebras have been implemented in Biology to represent different aspects of non-Mendelian inheritance (particularly, uniparental inheritance) in an algebraic way. Thus, Tian himself \cite{Tian2011} made use of these algebras to analyze how the homoplasmy of a cell population (represented by persistent generators) can derive from an heteroplasmy one (represented by transient generators). This is useful, for instance, to study mitochondrial disorders and mutations in tissues of patients. He also implemented evolution algebras to describe genetically dynamical patterns that provide information about the asexual reproduction of Phytophthora infectans causing the late blight of tomatoes and potatoes. An algebraic concept as nilpotency represents the extinction of original genetic types after certain generations.
	
	Continuous-time dynamical systems of mosquito populations were studied by Rozikov and Velasco \cite{RV}. Mosquito populations were already been dealt with from an algebraic point of view in 2011 by Junliang Lu and Jia Li \cite{Luli}, who proposed a discrete time structured model of such a population. On their own, Rozikov and Velasco considered a discrete-time dynamical system described by an evolution algebra with two fixed points, which become saddle points under some conditions on the parameters of the system. Its structure matrix is equal to the Jacobian of the QSO at a fixed point.
	
	Another aspect that Tian \cite{Tian2008} realized was the possible use of evolution algebras on {\em coalescent theory} \cite{Tian2005}. That is, on genetic evolution reversely over time.  The study of backwards evolution of Mendelian genetic systems by means of coalgebras had been previously introduced in 2004 by Tian and Bai-Lian Li \cite{TianLi}. In 2019, Paniello \cite{Paniello} introduced the concept of {\em evolution coalgebra} in order to model backwards evolution of Non-Mendelian genetic systems. She established the connection between such coalgebras and evolution algebras, by focusing in particular on the cases of genetic realizations. An illustrative example of this approach was given by Paniello herself \cite{Paniello2020a} in case of dealing with chicken populations.
	
	In 2020, Miguel Bustamante, Mellon and Velasco \cite{Bustamante2020, Bustamante2020a} described necessary and sufficient conditions for a given algebra to be an evolution algebra. They proved that this problem is equivalent to the simultaneous diagonalization via congruence of a given set of matrices. Based on this fact, they realized that arbitrarily small perturbations of classical genetic algebras representing Mendelian and auto-tetraploid inheritance (which are not evolution algebras) may give rise to evolution algebras. The algebraic representations of sexual and asexual inheritance by means of genetic and evolution algebras seem, therefore, to be closer than previously thought. A further study of baric and evolution algebras was established by the authors as a first stage to better understand this relationship. In any case, notice that Mukhamedov and Qaralleh already set out the following question in 2014. 
	
	\begin{prob}[\cite{Mukhamedov2014}] Is there a transformation of a given genetic algebra to some evolution algebra?
	\end{prob}
	
	They gave an affirmative answer in case of dealing with two-dimensional genetic algebras (see \cite[Theorem 4.1]{Mukhamedov2014}). The resulting transformation enabled them to provide necessary conditions on the structure matrix of these genetic algebras for ensuring the existence of non-trivial derivations (see \cite[Theorem 5.1]{Mukhamedov2014}).
	
	\section{Conclusion}
	
	This paper has dealt with the past, origin, development and possible further work of the theory of evolution algebras introduced by Tian \cite{tesistian} in 2004, with particular emphasis in the relationship that this theory has with a wide amount of distinct branches, not only in Mathematics, but with other areas of research. As such, the paper aims to be a starting point for all those researchers interested in this topic.
	
	\section*{Acknowledgements}
	
	The authors want to express their gratitude to the anonymous referee for the comprehensive reading of the paper and his/her pertinent comments and suggestions, which helped improve the manuscript.
	
	Falc\'on's work is partially supported by the research project FQM-016 from Junta de Andaluc\'\i a.

\end{document}